\documentclass[12pt]{article}
\usepackage{amsfonts}
\usepackage[dvips]{graphicx}

\begin{document}

\title{\bf An estimate of accuracy for interpolant numerical solutions of a PDE problem}

\date{December 2004}

\author{\bf Gianluca Argentini \\
\normalsize gianluca.argentini@riellogroup.com \\
\textit{Advanced Computing Laboratory}\\
Information \& Communication Technology Department\\
\textit{Riello Group}, 37045 Legnago (Verona), Italy}

\maketitle

\begin{abstract}
In this paper we present an estimate of accuracy for a piecewise polynomial approximation of a classical numerical solution to a non linear differential problem. We suppose the numerical solution $u$ is computed using a grid with a small linear step and interval time $T_u$, while the polynomial approximation $v$ is an interpolation of the values of a numerical solution on a less fine grid and interval time $T_v << T_u$. The estimate shows that the interpolant solution $v$ can be, under suitable hypotheses, a good approximation and in general its computational cost is much lower of the cost of the fine numerical solution. We present two possible applications to linear case and periodic case.
\end{abstract}

\section{The problem}

Let $[a,b]$, $a < b$, an interval on the real line $\mathbb{R}$, and let $C^1([0,+\infty)\times[a,b])$ the space of
continuously differentiable real functions defined on $[0,+\infty)\times[a,b]$. 
Let $F: C^1([0,+\infty)\times[a,b]) \rightarrow \mathbb{R}$ a continuously differentiable real functional.
Then, if $u\in C^1([0,+\infty)\times[a,b])$, the partial differential equation, usually associated to hyperbolic conservation laws,

\begin{equation}\label{pde1}
	\frac{\partial u}{\partial t} + \frac{\partial}{\partial x}F(u) = 0
\end{equation}

\noindent is in general non linear on the unknown $u$. Let $u_0 \in C^1([a,b])$ and $u_a$, $u_b$ two real numbers. We consider the
boundary problem

\begin{equation}
\left\{
\begin{array}{ll}\label{bvp1}
	$$\frac{\partial u}{\partial t}(x) + \frac{\partial F(u)}{\partial x}(x) = 0 \hspace{0.2cm} \forall x\in(a,b)$$\\
	$$u(t,a) = u_a$$, \hspace{0.1cm} $$u(t,b) = u_b \hspace{0.2cm} \forall t\in[0,+\infty)$$\\
	$$u(0,x) = u_0(x) \hspace{0.2cm} \forall x\in(a,b)$$
\end{array}
\right.
\end{equation}

\noindent and its numerical resolution. In this paper we consider only the {\it Euler forward in time - centered in space} ({\it EFC}) schema (see \cite{quva}) for Finite Differences method, but for other numerical schemas one can apply the same considerations with analogous, but more complicated, demonstrations. Let $h > 0$ a linear step for a partition $(x_j)_j$, $0\leq j \leq P$, $P=\frac{b-a}{h}$, of the interval $[a,b]$ such that $x_0=a$, $x_P=b$, and let $\Delta t$ the time step of the numerical schema. We denote by $u^n = (u_j^n)_j$ the numerical solution of problem (\ref{bvp1}) computed at time step $n$-th, hence at instant $t=n\Delta t$.\\

In {\it EFC} schema the values at time step $n+1$ are so computed (see \cite{an},{\cite{quva}):

\begin{equation}\label{eac1}
	u_j^{n+1} = u_j^n - F'(u_j^n)\frac{\Delta t}{2h}(u_{j+1}^n - u_{j-1}^n), \hspace{0.2cm} 0 < j < P-1
\end{equation}

\noindent where $F'(u_j^n)$ is the value of $\frac{dF}{du}(u)$ obtained by substitution of $u$ with the discretized value $u_j^n$. We suppose that the schema (\ref{eac1}) satisfies the {\it Courant-Friedrichs-Lewy condition} ({\it CFL}) (see \cite{quva}):

\begin{equation}\label{cfl1}
	\left|F'(u_j^n)\right|\frac{\Delta t}{h}\leq 1 \hspace{0.2cm} \forall j\in {\mathbb Z}, \forall n\in {\mathbb N}, \forall h\in {\mathbb R_+}
\end{equation}

\noindent This condition is usually admitted in the literature about numerical resolution of PDEs and is related to the question of {\it stability} of solutions.\\

Let $N\in {\mathbb N}$, $N>1$, $r\in {\mathbb N_+}$, $k=\frac{h}{r}$, $dt=\frac{\Delta t}{r}$, $u^N$ the numerical solution of (\ref{bvp1}) at time step $N$ with linear step $h$, $u^M$ the solution at time step $M=Nr$ with linear step $k$. Note that {\it CFL} condition is satisfied in the second case too, because $\frac{\Delta t}{h} = \frac{dt}{k}$. The solution $u^M$ is computed on a grid finer than that of $u^N$, being $h > k$, and hence it can be considered a better approximation of the exact solution, if exists, of the problem (\ref{bvp1}). But $u^M$ has the disadvantage of a much greater computational cost, which can be very high for $r >> 1$ (see \cite{ar}).\\

We construct a piecewise polynomial interpolation $v$ of $u^N$ on the nodes of the partition with grid step $h$. Function $v$ can be considered as a {\it global continuous and differentiable} representation of the discrete values $u^N$ on the interval $[a,b]$. Then one can compute, for every couple $(x_j,x_{j+1})$ of nodes, the interpolating polynomial $v_j = v_j(t)$ on a set of $r$ values of $t$. In the following, we present an estimate of the error $|v_j(t) - u_x^M|$, where $t=\frac{x-x_j}{h}$ (see Section 3), with the aim of establishing if $v$ can be considered a good approximation for the solution of the problem (\ref{bvp1}).

\section{Construction of interpolating polynomials}

Let $(x_j,x_{j+1})$ a couple of nodes internal to $[a,b]$, hence $1\leq j\leq P-2$. For simplicity we don't consider the boundary nodes $x_0$ and $x_P$, but with a suitable choice of a discretization formula for the spatial derivative $\frac{\partial u}{\partial x}$ (see \cite{an}, \cite{st}) we could apply analogous considerations. We write, for more simple notation, $p_1=x_j$, $p_2=x_{j+1}$, $d_1=u_j^N$, $d_2=u_{j+1}^N$. We construct a cubic polynomial $q=q(t)=at^3+bt^2+ct+d$, $0\leq t\leq 1$, imposing the following four conditions:

\begin{description}
	\item[] {\it 1}. \hspace{0.3cm} $q(0)=p_1$;
	\item[] {\it 2}. \hspace{0.3cm} $q(1)=p_2$;
	\item[] {\it 3}. \hspace{0.3cm} $q'(0)=d_1$; 
	\item[] {\it 4}. \hspace{0.3cm} $q'(1)=d_2$;
\end{description}

\noindent where $q'(t)=3at^2+2bt+c$. Then, solving the following linear system for the unknowns $a$, $b$, $c$, $d$

\begin{equation}
\left\{
\begin{array}{ll}\label{system1}
	d = p_1 \\
	a+b+c+d = p_2 \\
	c = d_1 \\
	3a+2b+c = d_2
\end{array}
\right.
\end{equation}

\noindent one obtains the expression for the derivative of each cubic:

\begin{equation}\label{cubic1}
	q'(t)=(6p_1-6p_2+3d_1+3d_2)t^2+(-6p_1+6p_2-4d_1-2d_2)t+d_1
\end{equation}

\section{Auxiliary technical results}

Let $w=u^N$ the numerical solution of problem (\ref{bvp1}) obtained with a $EFC$ schema with grid step $h$ ($h$-grid) and initial condition $u_0$ sampled using step $h$. We consider a couple of internal nodes $(x_j, x_{j+1})$ on $h$-grid, and a subpartition $(x_{j,i})_{0\leq i \leq r}$, where $r\in \mathbb{N}_+$, $mod(r,2)=0$, $x_j=x_{j,0},...,x_{j+1}=x_{j,r}$. Let $v=v_j(t)=q_j'(t)$, where $q_j$ is the cubic computed solving the system (\ref{system1}). We call $v$ an {\it interpolant solution} of the problem (\ref{bvp1}). We define the set $t_m = \frac{m}{r}$, where $m\in \mathbb{N}$, $0\leq m \leq r$. If $s=\frac{x_{j+1}-x_j}{r}$, we have $x_m=x_{j,m}=ms+x_j$, therefore

\begin{equation}
	t_m = \frac{m}{r} = \frac{m}{r}\frac{x_{j+1}-x_j}{h} = m\frac{x_{j+1}-x_j}{r}\frac{1}{h} = \frac{ms}{h}=\frac{x_m-x_j}{h}
\end{equation}

Let $k=\frac{h}{r}$, and let $u=u^M$ the solution obtained with grid step $k$ and $u_0$ sampled using $k$. For simplicity we denote $u_m$ the value of $u$ on the node $x_m = ms + x_j$, where $j$, as before, is index of the $h$-grid. In the following we demonstrate some preliminary results for estimate the error $|v_m - u_m|$, where $v_m=v(t_m)$.\\

\newtheorem{prop}{Proposition}

\begin{prop}\label{prop1}
	If $t\in [0,1]$, exists a value $s\in [0,1]$ such that
	\begin{equation}\label{s1}
		t^2 = t - \frac{s}{4}
	\end{equation}
\end{prop}

{\it Proof}. If $y=t^2-t$, then $y\leq 0$ and $y'=2t-1$, therefore $t_0=\frac{1}{2}$ is the solution of $y'(t)=0$. At $t_0$ the function $y$ has a minimum, because $y''=2>0$. The value is $y(\frac{1}{2})=-\frac{1}{4}$. Therefore we have $-\frac{1}{4}\leq t^2-t\leq 0$, from which (\ref{s1}) follows with $s\in [0,1]$.\hspace{0.2cm} $\square$\\

\begin{prop}\label{prop2}
	For every $m$ such that $0\leq m \leq r$:
	\begin{equation}\label{v1}
		|v_m - v_r|\leq \frac{3}{2}h + 2|d_1-d_2| + \frac{3}{4}|d_1+d_2|
	\end{equation}
\end{prop}

{\it Proof}. From (\ref{cubic1}) we have $v_r = v(1) = d_2$, therefore

\begin{eqnarray}
\nonumber{ |v_m - v_r|\leq |(6p_1-6p_2+3d_1+3d_2)t^2+(-6p_1+6p_2-4d_1-2d_2)t+d_1-d_2| } \\
\nonumber{ \leq |(6p_1-6p_2)(t^2-t)| + |(3d_1+3d_2)t^2 - (4d_1+2d_2)t| + |d_1-d_2| }
\end{eqnarray}

\noindent Then, using (\ref{s1}), we obtain

\begin{equation}
\nonumber{ |v_m - v_r|\leq \frac{3}{2}|p_1-p_2| + 2|d_1-d_2| + \frac{3}{4}|d_1+d_2| }
\end{equation}

\noindent and hence (\ref{v1}), because $|p_1-p_2|=h$.\hspace{0.2cm} $\square$\\

In the same manner, using the fact that $v_0 = v(0) = d_1$, one proofs the following\\

\begin{prop}\label{prop3}
	For every $m$ such that $0\leq m \leq r$:
	\begin{equation}\label{v2}
		|v_m - v_0|\leq \frac{3}{2}h + 2|d_1-d_2| + \frac{3}{4}|d_1+d_2|
	\end{equation}
\end{prop}

In the following propositions $\epsilon$ is a fixed positive real number, $w^0$ is the initial condition $u_0$ of (\ref{bvp1}) sampled on the nodes of $h$-grid and $u^0$ the initial condition sampled on $k$-grid nodes.\\

\begin{prop}\label{prop4}
	If $|u_i^0-u_{i+1}^0|\leq\frac{\epsilon}{3^M}$ \hspace{0.05cm} $\forall i$, then $|u_i^n-u_{i+1}^n|\leq\frac{\epsilon}{3^{(M-n)}}$ \hspace{0.05cm} $\forall i$, \hspace{0.1cm} $\forall n\leq M, n\in \mathbb{N}_+$.
\end{prop}

{\it Proof}. By induction. Let $n=1$. We have

\begin{equation}
	|u_i^1-u_{i+1}^1| \leq |u_i^1-u_i^0|+|u_i^0-u_{i+1}^0|+|u_{i+1}^0-u_{i+1}^1|
\end{equation}

\noindent Using the {\it EFC} schema (\ref{eac1}), the {\it CFL} condition (\ref{cfl1}) and the hypothesis, we have
\begin{eqnarray}
	\nonumber{ |u_i^1-u_{i+1}^1| \leq \frac{1}{2}|u_{i+1}^0-u_{i-1}^0|+|u_i^0-u_{i+1}^0|+\frac{1}{2}|u_{i+2}^0-u_i^0| }\\
	\nonumber{ \leq \frac{1}{2}|u_{i+1}^0-u_i^0|+\frac{1}{2}|u_i^0-u_{i-1}^0|+|u_i^0-u_{i+1}^0|+\frac{1}{2}|u_{i+2}^0-u_{i+1}^0|+\frac{1}{2}|u_{i+1}^0-u_i^0| }\\
	\nonumber{ \leq 3\frac{\epsilon}{3^M} = \frac{\epsilon}{3^{(M-1)}} }
\end{eqnarray}

\noindent that is the first inductive step. The intermediate inductive  step is $|u_i^n-u_{i+1}^n|\leq\frac{\epsilon}{3^{(M-n)}}$. Then
\begin{eqnarray}
	\nonumber{ |u_i^{n+1}-u_{i+1}^{n+1}| \leq |u_i^{n+1}-u_i^n|+|u_i^n-u_{i+1}^n|+|u_{i+1}^n-u_{i+1}^{n+1}| }\\
	\nonumber{ \leq \frac{1}{2}|u_{i+1}^n-u_{i-1}^n|+|u_i^n-u_{i+1}^n|+\frac{1}{2}|u_{i+2}^n-u_i^n| }\\
	\nonumber{ \leq \frac{1}{2}|u_{i+1}^n-u_i^n|+\frac{1}{2}|u_i^n-u_{i-1}^n|+|u_i^n-u_{i+1}^n|+\frac{1}{2}|u_{i+2}^n-u_{i+1}^n|+\frac{1}{2}|u_{i+1}^n-u_i^n| }\\
	\nonumber { \leq 3\frac{\epsilon}{3^{(M-n)}} = \frac{\epsilon}{3^{(M-(n+1))}} }
\end{eqnarray}

\noindent and this is the final inductive step. \hspace{0.2cm}$\square$\\

\begin{prop}\label{prop5}
	If $|u_i^0-u_{i+1}^0|\leq\frac{\epsilon}{3^M}$ \hspace{0.05cm} $\forall i$, then $|u_i^{n+1}-u_i^n|\leq\frac{\epsilon}{3^{(M-n)}}$ \hspace{0.05cm} $\forall i$, \hspace{0.05cm} $\forall n\leq M$.
\end{prop}

{\it Proof}. Let $n\geq 0$. Then
\begin{eqnarray}
	\nonumber{ |u_i^{n+1}-u_i^n|\leq \frac{1}{2}|u_{i+1}^n-u_{i-1}^n| } \\
	\nonumber{ \leq \frac{1}{2}(|u_{i+1}^n-u_i^n|+|u_i^n-u_{i-1}^n|) }
\end{eqnarray}

\noindent and hence, from Proposition (\ref{prop4}) 

\begin{equation}
	|u_i^{n+1}-u_i^n|\leq \frac{\epsilon}{3^{(M-n)}}
\end{equation}

\noindent that is the thesis. \hspace{0.2cm}$\square$\\

In the next proposition, the index $i$ is associated to the node $x_i$ in the $k$-grid, while the index $j$ is associated to the node $x_j$ in the $h$-grid.\\

\begin{prop}\label{prop6}
	If $|u_i^0-u_{i+1}^0|\leq\frac{\epsilon}{3^M}$ \hspace{0.05cm} $\forall i$, then $|u_j^0-u_{j+1}^0|\leq\frac{\epsilon}{3^N}$ \hspace{0.05cm} $\forall j$.
\end{prop}

{\it Proof}. Using the fact that $r < 3^r$ and $M=Nr \geq (N+r)$ if $r,N > 1$, we have
\begin{eqnarray}
	\nonumber{ |u_j^0-u_{j+1}^0| \leq \sum_{m=0}^{r-1}|u_{j+\frac{m}{r}}^0-u_{j+\frac{m+1}{r}}^0| \leq \sum_{m=0}^{r-1} \frac{\epsilon}{3^{Nr}} \leq \sum_{m=0}^{r-1} \frac{\epsilon}{r3^N} = \frac{\epsilon}{3^N} }
\end{eqnarray}

\noindent therefore the proposition is verified. \hspace{0.2cm}$\square$\\

Note that by definition $w_j^0=u_j^0$, hence from Proposition (\ref{prop6}) follows that

\begin{equation}\label{w0ineq1}
	|w_j^0-w_{j+1}^0|\leq\frac{\epsilon}{3^N}
\end{equation}

Then the Propositions (\ref{prop4}) and (\ref{prop5}) hold for the solutions ($w^k$) too, using $N$ instead of $M$:\\

\newtheorem{corol}{Corollary}

\begin{corol}\label{corol1}
	If $|u_i^0-u_{i+1}^0|\leq\frac{\epsilon}{3^M}$ \hspace{0.05cm} $\forall i$, then
  $|w_j^{k+1}-w_j^k|\leq\frac{\epsilon}{3^{(N-k)}}$, $|w_j^k-w_{j+1}^k|\leq\frac{\epsilon}{3^{(N-k)}}$ \hspace{0.05cm} $\forall j$, \hspace{0.05cm} $\forall k\leq N$.
\end{corol}

In particular, for $k=N$, writing $w_j^N=d_1$ and $w_{j+1}^N=d_2$, follows that\\

\begin{corol}\label{corol2}
	If $|u_i^0-u_{i+1}^0|\leq\frac{\epsilon}{3^M}$ \hspace{0.05cm} $\forall i$, then
  $|d_1-d_2|\leq \epsilon$.
\end{corol}

In the next proposition the simplified symbol $u_0 = u_j^M$ is the value of solution $u$ at node $x_{jr}$ of the $k$-grid, i.e. the node $x_j$ of the $h$-grid.\\

\begin{prop}\label{prop7}
	If $|u_i^0-u_{i+1}^0|\leq\frac{\epsilon}{3^M}$ \hspace{0.05cm} $\forall i$, then $|v_0-u_0|\leq \epsilon$.
\end{prop}

{\it Proof}. By definition, $v_0 = d_1 = w_j = w_j^N$, where $w^N$ is the solution at time step $N$ of (\ref{bvp1}) with initial condition $w_0$. Then, using a telescopic sum like $(a_N-a_0)=\sum_{k=1}^N(a_k-a_{k-1})$:

\begin{equation}\label{prop7.1}
	|v_0-u_0|=|w_j-u_0| \leq \sum_{k=1}^N |w_j^k-w_j^{k-1}| + |w_j^0-u_j^0| + \sum_{h=1}^M |u_j^{h-1}-u_j^h|
\end{equation}

\noindent By definition, $w_j^0=u_j^0$, hence the second term of the second member is null. For the first sum, from hypothesis and from Corollary (\ref{corol1}) follows
\begin{eqnarray}
	\nonumber{ \sum_{k=1}^N |w_j^k-w_j^{k-1}| \leq \sum_{k=1}^N \frac{\epsilon}{3^{(N-k+1)}} } \\
	\nonumber{ = \frac{\epsilon}{3^{(N+1)}} \sum_{k=1}^N 3^k = \frac{\epsilon}{3^{(N+1)}} \frac{3^{(N+1)}-1}{2} \leq \frac{\epsilon}{2}  }
\end{eqnarray}

\noindent For the second sum the computation, using the Proposition (\ref{prop5}), is the same. Therefore the proposition is verified. \hspace{0.2cm}$\square$\\

Using the simplified symbol $u_r = u_{(j+1)r}^M$ for the value of $u$ at the node $x_{(j+1)r}$ at time step $M$, in the same manner we can proof the following\\

\begin{prop}\label{prop8}
	If $|u_i^0-u_{i+1}^0|\leq\frac{\epsilon}{3^M}$ \hspace{0.05cm} $\forall i$, then $|v_r-u_r|\leq \epsilon$.\\
\end{prop}

In the next Proposition we remember that $u_m$ is the value of $u=u^M$ on the node $x_m = ms + x_j$.\\

\begin{prop}\label{prop9}
	Let $1\leq m\leq \frac{r}{2}$. If $|u_i^0-u_{i+1}^0|\leq\frac{\epsilon}{3^M}$ \hspace{0.05cm} $\forall i$, then $|u_0-u_m|\leq m\epsilon$.\\
\end{prop}

{\it Proof.} Using Proposition (\ref{prop4}) with $n=M$, we have
\begin{eqnarray}
	\nonumber{ |u_0-u_m| \leq \sum_{k=1}^m |u_{k-1}^M-u_k^M| \leq m{\epsilon} \hspace{0.2cm}\square}
\end{eqnarray}

\begin{prop}\label{prop10}
	Let $\frac{r}{2}\leq m\leq r$. If $|u_i^0-u_{i+1}^0|\leq\frac{\epsilon}{3^M}$ \hspace{0.05cm} $\forall i$, then $|u_m-u_r|\leq (r-m)\epsilon$.\\
\end{prop}

{\it Proof.} Using Proposition (\ref{prop4}) with $n=M$, we have
\begin{eqnarray}
	\nonumber{ |u_m-u_r| \leq \sum_{k=1}^{r-m} |u_{m+k-1}^M-u_{m+k}^M| \leq (r-m){\epsilon} \hspace{0.2cm}\square}
\end{eqnarray}

\section{The main results}

Now we can stated the main results on the estimate about the accuracy of the interpolant solution $v$.

\newtheorem{theo}{Theorem}

\begin{theo}\label{theo1}
Let $\epsilon > 0$, $m\in \mathbb{N}_+$, $1\leq m\leq r$. If $|u_i^0-u_{i+1}^0|\leq\frac{\epsilon}{3^M}$ \hspace{0.05cm} $\forall i$, then
\begin{equation} 
	|v_m-u_m| \leq \frac{3}{2}h + \frac{3}{4}|d_1+d_2| + (min[m,r-m]+3)\epsilon
\end{equation}
\end{theo}
\vspace{0.5cm}

{\it Proof}. If $m \leq \frac{r}{2}$, let $p=0$, else $p=r$. Using Propositions (\ref{prop2}) or (\ref{prop3}), (\ref{prop7}) or (\ref{prop8}), (\ref{prop9}) or (\ref{prop10}) and Corollary (\ref{corol1}), we have
\begin{eqnarray}
	\nonumber{ |v_m-u_m| \leq |v_m-v_p| + |v_p-u_p| + |u_p-u_m| }\\
	\nonumber{ \leq \frac{3}{2}h + 2\epsilon + \frac{3}{4}|d_1+d_2| + (min[m,r-m]+1)\epsilon }
\end{eqnarray}

\noindent because if $m\leq \frac{r}{2}$ then $min[m,r-m]=m$, else $min[m,r-m]=r-m$.$\hspace{0.2cm}$ $\square$ \\

\begin{corol}\label{corol3}
Let $\epsilon > 0$, $m\in \mathbb{N}_+$, $1\leq m\leq r$. There exists $\delta > 0$ such that if $|u_i^0-u_{i+1}^0|\leq\delta$ \hspace{0.05cm} $\forall i$, then
\begin{equation}\label{corol3dis}
	|v_m-u_m| \leq \frac{3}{2}h + \frac{3}{4}|d_1+d_2| + \epsilon
\end{equation}
\end{corol}
\vspace{0.5cm}

{\it Proof}. Let $\mu = \frac{\epsilon}{min[m,r-m]+3}$. Then $\mu > 0$, and from Theorem (\ref{theo1}), defined $\delta = \frac{\mu}{3^M}$, follows that if $|u_i^0-u_{i+1}^0|\leq\delta$, then

\begin{equation}
	|v_m-u_m| \leq \frac{3}{2}h + \frac{3}{4}|d_1+d_2| + (min[m,r-m]+3)\mu
\end{equation}

\noindent But $(min[m,r-m]+3)\mu = \epsilon$, so the thesis is verified.$\hspace{0.2cm}\square$ \\

The estimate stated in Corollary (\ref{corol3}) is of the same kind of the estimates established, in the case of finite element method, by Kuznetsov and Tadmor-Tang for approximated numerical solutions to scalar conservation laws equations (see \cite{tata}).\\
In particular, inequality (\ref{corol3dis}) shows that the accuracy, respect to the fine $u$ solution computed using a grid step $k$, of the interpolant solution $v$, computed using a grid step $h=rk$, is of the first order on $h$. In (\ref{corol3dis}) the role of the time steps $N$ and $M$ is hidden in the third addendum $\epsilon$, as stated by Corollary (\ref{corol2}).\\
The next theorem gives a sufficient condition on the initial function $u_0$ for the validity of the hypothesis in Theorem (\ref{theo1}).\\

\begin{theo}\label{theo2}
Let $\epsilon > 0$, $h > 0$, $N\in \mathbb{N}$, $N>1$. If $u_0 \in C^0([a,b])$, then exists $r\in \mathbb{R}_+$ such that, if $0\leq m\leq r$,
\begin{equation} 
	|v_m-u_m| \leq \frac{3}{2}h + \frac{3}{4}|d_1+d_2| + \epsilon
\end{equation}
\end{theo}
\vspace{0.5cm}

{\it Proof}. The function $u_0$, continuous on a closed and limited subset of $\mathbb{R}$, is uniformly continuous, therefore exists $\delta > 0$ such that if $x_1, x_2 \in[a,b]$, $|x_1-x_2|\leq \delta$, then $|u_0(x_1)-u_0(x_2)|\leq \epsilon$. We can choose the parameter $r\in \mathbb{N}_+$, and hence the time step $M=rN$, such that $\frac{\epsilon}{3^M}\leq \delta$. Therefore, with this choice of $r$, it can be applied the theorem (\ref{theo1}) and, with the right choice of $\delta$, the Corollary (\ref{corol3}).$\hspace{0.2cm}\square$ \\

In the next section we discuss the role of the term $|d_1+d_2|$ in some particular cases.

\section{Applications to particular cases}

As first application of the previous estimate, let $F(u)=au$, with $a\in \mathbb{R}$, $a\neq0$. Hence we suppose that $F$ is linear.
In this case the numerical schema EFC is {\it stable} in the norm $|u|_2=\left(h \sum_j|u_j|^2 \right)^{\frac{1}{2}}$ until time step $N$, i.e. $|u|_2 \leq C_N|u_0|_2$ with $C_N=C(N)\in \mathbb{R}_+$, if $\Delta t \leq \left(\frac{h}{a}\right)^2$; the constant $C_N$ is equal to $e^{\frac{N\Delta t}{2}}$ (see \cite{quva}), while $j$ is the index on $h$-grid. If we choose the grid step $h$ such that $h\leq |a|$, than the $CFL$ condition is satisfied because

\begin{equation}
	|a|\frac{\Delta t}{h}\leq \frac{h}{|a|} \leq 1
\end{equation}

\newtheorem{lemma}{Lemma}

\begin{lemma}\label{lemma1}
 Let $m\in \mathbb{N}$, $m\geq 1$, and $(C_j)_{0\leq j\leq m}$ a set of positive real numbers. If $A$, $B$ $\in (C_j)_j$, then
 
 \begin{equation}
 	A+B\leq 2\left(\sum_{j=0}^m C_j^2\right)^{\frac{1}{2}}
 \end{equation}

\end{lemma}

{\it Proof}. We have $A^2\leq \sum_{j=0}^m C_j^2$ and $B^2\leq \sum_{j=0}^m C_j^2$, therefore

\begin{equation}
 A\leq \left(\sum_{j=0}^m C_j^2\right)^{\frac{1}{2}}, \hspace{0.5cm} B\leq \left(\sum_{j=0}^m C_j^2\right)^{\frac{1}{2}}
\end{equation}

\noindent hence the thesis is verified.$\hspace{0.2cm}\square$\\

If $u \in C^0([a,b])$, let $|u|_\infty=sup_{x\in[a,b]}|u(x)|$.

\begin{corol}\label{corol4}
	Let $\epsilon > 0$, $N\in \mathbb{N}$, $N>1$, $h > 0$, $a\neq0$, $h\leq |a|$, $\Delta t \leq \left(\frac{h}{a}\right)^2$ and $u_0 \in C^0([a,b])$. Then exists $r\in \mathbb{R}_+$ such that, if $0\leq m\leq r$,
\begin{equation} 
	|v_m-u_m| \leq \frac{3}{2}\left[h + \left(\frac{b-a}{h}+1\right)^{\frac{1}{2}}e^{\frac{N\Delta t}{2}}|u_0|_\infty\right] + \epsilon
\end{equation}
\end{corol}
\vspace{0.5cm}

{\it Proof}. From theorem (\ref{theo2}) we have

\begin{equation} 
	|v_m-u_m| \leq \frac{3}{2}h + \frac{3}{4}|d_1+d_2| + \epsilon
\end{equation}

\noindent where $d_1=w_j^N$ and $d_2=w_{j+1}^N$, for a generic node $x_j$ of $h$-grid. Let $s=\frac{b-a}{h}$ and $|w_0|_\infty = sup_j|w_{0,j}|$. From Lemma (\ref{lemma1}) and from stability in the norm $|\cdot|_2$ follows that

\begin{eqnarray}
	\nonumber{ \frac{3}{4}|d_1+d_2| \leq \frac{3}{4}(|w_j^N|+|w_{j+1}^N|) \leq \frac{3}{2} \left(\sum_{j=0}^s (w_j^N)^2\right)^{\frac{1}{2}} }\\
	\nonumber{ = \frac{3}{2h^{\frac{1}{2}}}|w^N|_2 \leq \frac{3}{2h^{\frac{1}{2}}}e^{\frac{N\Delta t}{2}}|w_0|_2 }\\
	\nonumber{ \leq \frac{3}{2}\left(s+1\right)^{\frac{1}{2}}e^{\frac{N\Delta t}{2}}|w_0|_\infty }
\end{eqnarray}

\noindent and the thesis is verified because $|w_0|_\infty\leq|u_0|_\infty$.$\hspace{0.2cm}\square$ \\

The Corollary (\ref{corol4}) shows that the error $|v_m-u_m|$ is {\it uniformly limited} on interval $[a,b]$ for every time step $n\leq N$. The presence of some kind of norm on $u_0$, like $|u_0|_\infty$, in the upper limitation is typical in literature about accuracy of numerical solutions (see \cite{tata}).\\

As second application, let $w$ the numerical solution at time step $N$ of the problem (\ref{bvp1}) on the $h$-grid, and let $x_j$, $x_{j+1}$ two nodes of this grid. We call the interval $[x_j,x_{j+1}]$ an {\it interval of rising turbulence at time step $N$} for the $h$-grid if $(w_jw_{j+1}) < 0$. The definition is based on the intuitive fact that a change of sign for the velocity field of a linear flow induces a superposition of particles in the motion of a fluid. In one dimensional case, in this situation the velocity $u$ in the problem (\ref{bvp1}) might be not a classic function.\\

\begin{corol}\label{corol5}
	Let $\epsilon > 0$, $N\in \mathbb{N}$, $N>1$, $h > 0$, $u_0 \in C^0([a,b])$ and $[x_j,x_{j+1}]\subset(a,b)$ an interval of rising turbulence at time step $N$ for the $h$-grid. Let $w^N$ the numerical solution at time step $N$. Then exists $r\in \mathbb{R}_+$ such that, if $0\leq m\leq r$, $v_0=w_j^N$ and $v_r=w_{j+1}^N$,
\begin{equation} 
	|v_m-u_m| \leq \frac{1}{2}(3h + 5\epsilon)
\end{equation}
\noindent where $u$ is the solution on $k$-grid, $k=rh$.
\end{corol}
\vspace{0.5cm}

{\it Proof}. From theorem (\ref{theo2}) we have

\begin{equation}\label{corol5eqn}
	|v_m-u_m| \leq \frac{3}{2}h + \frac{3}{4}|w_j^N+w_{j+1}^N| + \epsilon
\end{equation}

\noindent From hypothesis we suppose $w_j=w_j^N<0$ and $w_{j+1}=w_{j+1}^N>0$. Then exists $w>0$ such that $w_j=-ww_{j+1}$. Theorem (\ref{theo2}) holds, and from Corollary (\ref{corol1}) we have $|w_j-w_{j+1}|\leq \epsilon$. Therefore

\begin{equation}
	|w_j-w_{j+1}| = |-ww_{j+1}-w_{j+1}|=w_{j+1}|-w-1|=w_{j+1}(1+w) \leq \epsilon
\end{equation}

\noindent and hence $w_{j+1}\leq \frac{\epsilon}{1+w}\leq \epsilon$. Also

\begin{equation}
	|w_j-w_{j+1}| = -w_j+w_{j+1}\leq \epsilon
\end{equation}

\noindent hence $|w_j|=-w_j\leq \epsilon - w_{j+1}\leq \epsilon$. Therefore we have

\begin{equation}
	|w_j+w_{j+1}|\leq |w_j|+|w_{j+1}|\leq 2\epsilon
\end{equation}

\noindent and from (\ref{corol5eqn}) follows that

\begin{equation}
		|v_m-u_m| \leq \frac{3}{2}h + \frac{3}{2}\epsilon + \epsilon
\end{equation}

\noindent The thesis is verified.$\hspace{0.2cm}\square$ \\

From the proof of the Corollary we may note that in an interval of rising turbolence the values of the solution at extremes are small. The particular limit case $w_j=w_{j+1}=0$, for which, from theorem (\ref{theo2}), $|v_m-u_m| \leq \frac{3}{2}h + \epsilon$, might be interesting for searching periodic solution $u$ for which $u_j=u_{j+1}=0$ on the nodes of the $h$-grid:

\begin{figure}[ht]
	\begin{center}
	\includegraphics[width=13cm]{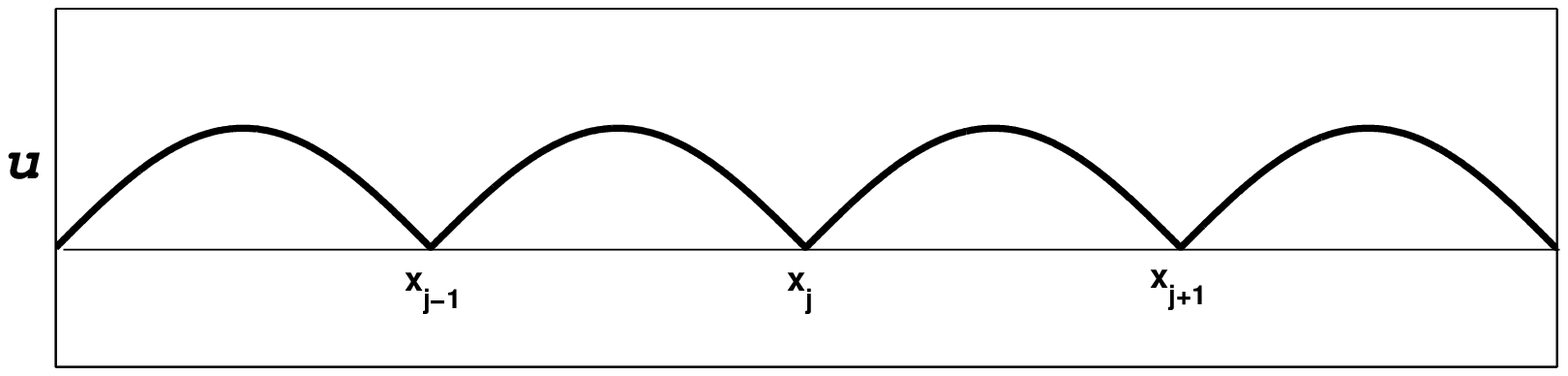}
	\end{center}
\end{figure}


\begin{thebibliography}{9}

\bibitem{an} J. D. Anderson, {\it Computational Fluid Dynamics},
McGraw-Hill, New York (1995).

\bibitem{ar} G. Argentini, {\it Using sparse matrices and splines-based interpolation in computational fluid dynamics simulation}, to be appeared in {\it Applied and Industrial Mathematics in Italy}, Acts of SIMAI 2004 Congress of Societ\`{a} Italiana di Matematica Applicata e Industriale, World Scientific Publishing Company, Singapore (2005).

\bibitem{qusasa} A. Quarteroni, R. Sacco and F. Saleri, {\it Numerical 
Mathematics}, volume {\bf 37} of {\it Texts in Applied Mathematics},
Springer-Verlag, New York (2000).

\bibitem{quva} A. Quarteroni and A. Valli, {\it Numerical Approximation of Partial
Differential Equations}, Springer, Berlin and Heidelberg (1994).

\bibitem{st} J. Strickwerda, {\it Finite Difference Schemes and Partial Differential Equations}, Wadworth \& Brooks/Cole, Pacific Grove, 1989.

\bibitem{tata} E. Tadmor and T. Tang, {\it Pointwise error estimates for relaxation approximations to conservations laws}, SIAM J. Math. Anal., {\bf 32}(4), pp. 870-886, (2000)

\end{thebibliography}
\end{document}